# Méthodologie de dimensionnement par optimisation bayésienne d'une machine synchro-réluctante assistée d'aimants permanents

Adán REYES REYES[1,2], André NASR[1], Delphine SINOQUET[1] and Sami HLIOUI[3]

[1]IFP Energies nouvelles, 1 et 4 avenue de Bois-Préau, 92852 Rueil-Malmaison, France ; Institut Carnot IFPEN Transports Energie.

[2]SATIE, CNRS, Université Paris-Saclay, ENS Paris-Saclay, Gif-sur-Yvette, France.

[3]SATIE, CNRS, Université Cergy Paris, Université Paris-Saclay, ENS Paris-Saclay, Gif-sur-Yvette, France.

**RESUME** - Dans cet article, trois approches d'optimisation différentes sont exploitées pour améliorer les performances d'une machine synchro-réluctante assistée par aimants permanents : une première optimisation utilisant des modèles de substitution fixes et deux approches d'optimisation bayésienne basées sur des modèles de substitution adaptatifs. Les résultats montrent que les approches bayésiennes aboutissent à des machines avec de meilleures performances, en utilisant le même temps de calcul (même nombre de simulations par éléments finis). Contrairement aux méthodologies d'optimisation basées sur des modèles de substitution fixes, les approches bayésiennes fournissent des solutions directement basées sur des simulations par éléments finis et ne nécessitent donc pas de vérifications a posteriori.

*Mots-clés* —*Machines synchrones, moteurs à aimant permanent, optimisation bayésienne, processus gaussiens, méthode des éléments finis.*

1. INTRODUCTION

L'optimisation est devenue de nos jours une étape primordiale dans le processus de dimensionnement des machines électriques. Cette étape consiste à trouver le bon compromis entre l'ensemble des objectifs visés, tout en respectant les différentes contraintes qui peuvent s'imposer. Dans les applications de mobilité électrique, on cherche souvent à maximiser le rendement du moteur électrique tout en minimisant le poids et le coût. Le dimensionnement par optimisation se présente alors comme une stratégie efficace et un outil de prise de décision.

Un problème d'optimisation de machines électriques est un problème complexe. D'une part, et à cause du comportement non-linéaire des matériaux magnétiques, le modèle électromagnétique utilisé pour simuler les performances de ces machines est aussi non linéaire. Cela conduit, par exemple, à l'utilisation de la méthode des Eléments Finis (EF) surtout que cette méthode présente aussi d'autres avantages comme sa précision et sa généricité (applicable sur n'importe quel dessin). D'autre part, les variables de conception impliquées dans cette optimisation sont nombreuses. Les concepteurs ont souvent recours à des algorithmes d'optimisation stochastiques tels que des algorithmes génétiques et l'optimisation par essaims particulaires [1][2] pour explorer l'espace de conception afin d'échapper aux optima locaux. Cependant, ces méthodes nécessitent un grand nombre de simulations EF, ce qui engendre un coût de calcul important.

Les approches classiquement utilisées pour réduire le temps d'une optimisation reposent sur des modèles de substitution des sorties d'intérêt du simulateur. Dans [3], des modèles de substitution fixes ont été construits à partir d'un échantillon limité de simulations EF et utilisés dans un algorithme d'optimisation multi-objectifs convergeant rapidement vers une approximation du front de Pareto « réel ». Cependant, des différences ont été constatées entre les résultats prédits par les modèles de substitution et ceux prédits par le simulateur EF. Afin d'améliorer la précision d'une telle optimisation, les modèles de substitution doivent être mis à jour à partir de nouvelles simulations pendant le processus d'optimisation, par exemple en utilisant une approche adaptative comme l'Optimisation Bayésienne (BO) [4][5].

Dans ce travail, nous avons appliqué trois méthodologies d'optimisation au problème de dimensionnement d'une machine électrique synchro-réluctante assistée d'aimants (PMaSynRel) : une première optimisation utilisant une approche basée sur des modèles de substitution fixes (Optim1), et deux approches d'optimisation bayésienne (Optim2 et Optim3). Les résultats issus de ces optimisations sont par la suite comparés et analysés. Dans la section suivante, nous introduirons la topologie de la machine ainsi que les paramètres de conception. Ensuite, nous allons décrire les processus d'optimisation et comparer les résultats obtenus.

2. GEOMETRIE DE LA MACHINE ETUDIEE ET VARIABLES D'OPTIMISATION

La machine étudiée dans cet article est illustrée à la Fig. 1. Il s'agit d'une machine triphasée synchro-réluctante assistée d'aimants à 8 pôles, 48 encoches et 3 aimants permanents par pôle. Le Tableau 1 présente les variables utilisées dans les optimisations. On peut noter que ces variables sont normalisées entre 0 et 1. Pour calculer la valeur réelle d'une variable, nous appliquons la formule suivante :

$$X_{réel} = X_{réel-min} + (X_{réel-max} - X_{réel-min}) * X_{norm} \quad (1)$$

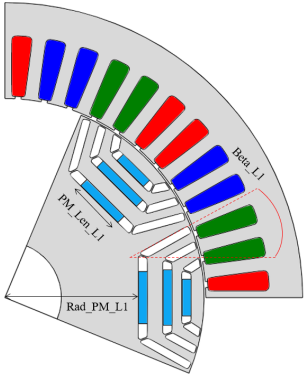

Fig. 1. Géométrie de la PMaSynRel.

Tableau 1 Variables d'optimisation

| Variables d'optimisation | Description | Borne Min Norm | Borne Max Norm |
|---|---|---|---|
| Rad_PM_L1 | Couche 1 : Rayon PM | 0.6 | 0.9 |
| Rad_PM_L3 | Couche 2 : Rayon PM | 0.1 | 0.9 |
| Rad_PM_L3 | Couche 3 : Rayon PM | 0.1 | 0.9 |
| Rad_Brid_L1 | Couche 1 : Pont radial | 0.1 | 0.9 |
| Rad_Brid_L2 | Couche 2 : Pont radial | 0.1 | 0.9 |
| Rad_Brid_L3 | Couche 3 : Pont radial | 0.1 | 0.9 |
| Beta_L1 | Couche 1: Angle d'ouverture | 0.1 | 0.9 |
| Beta_L2 | Couche 2: Angle d'ouverture | 0.1 | 0.9 |
| Beta_L3 | Couche 3: Angle d'ouverture | 0.1 | 0.9 |
| PM_Len_L1 | Couche 1: Longueur de l'aimant permanent | 0.1 | 0.9 |
| PM_Len_L2 | Couche 2: Longueur de l'aimant permanent | 0.1 | 0.9 |
| PM_Len_L3 | Couche 3: Longueur de l'aimant permanent | 0.1 | 0.9 |

Pour cette optimisation il existe trois bornes fixes : Beta_L1_min = 10°, Beta_L2_min = 10° et Beta_L3_min = 10°. A préciser que les paramètres géométriques du stator ont été fixés et seul le rotor a été optimisé.

3. OPTIMISATION D'UNE MACHINE PMASYNREL

Les trois optimisations seront utilisées pour maximiser le couple moyen de la machine ainsi que son rapport de puissance qui est défini comme le rapport entre la puissance à la vitesse maximale et la puissance maximale. Ce dernier objectif vise à obtenir des machines qui présentent une faible chute de puissance entre la vitesse de base et la vitesse maximale. Une seule contrainte limitant la FEM à vitesse maximale sera utilisée. Le problème d'optimisation multi-objectifs sous contraintes à résoudre s'écrit comme suit :

$$\min_{\mathbf{x}\in D_x} \{f_1(\mathbf{x}), f_2(\mathbf{x})\} \text{ tel que } g(\mathbf{x}) \leq 0 \quad (2)$$

où $f_1$ et $f_2$ sont les opposés du couple moyen et du rapport de puissance (afin de les maximiser), et g la valeur maximale de la FEM (peak) à la vitesse maximale (8000 RPM) moins un seuil de 650 $V_{pp}$. $D_x$ est l'espace des variables d'optimisation présenté dans le Tableau 1. Le Tableau 2 résume les objectifs et contraintes utilisés.

Tableau 2: Objectifs et contrainte de l''optimisation

| Fonction Objective 1 $f_1$ | Maximiser le couple moyen |
|---|---|
| Fonction Objective 1 $f_2$ | Minimiser les ondulations couple |
| Contrainte $g$ | FEM à vitesse maximal ≤ 650 $V_{pp}$ |

### 3.1. Construction des méta-modèles

Pour construire le métamodèle d'une fonction f, nous supposons que f est une réalisation d'un processus Gaussien entièrement défini par sa fonction moyenne, $\mu(x)$ et sa fonction de covariance $C(x, x')$ qui mesure la corrélation entre les valeurs de la fonction f entre deux points différents x et x'. Si $x = x'$, $C(x, x)$ est égale à la variance de densité de probabilité Gaussienne en x qui estime l'erreur de prédiction du métamodèle en ce point. Partant d'un plan d'expériences initial $D_n = \{x_i, f(x_i)\}_{i=1}^n$, le modèle de substitution $\hat{f}|D_n$ (noté $\hat{f}$ pour simplicité) est la moyenne $\mu_n(x)$ du GP conditionné aux données observées $D_n$, et sa variance donne une estimation des erreurs de prédiction. Plus de détails sur la construction des métamodèles du type GP pourront être trouvés dans [6].

Dans cette étude, la fonction de covariance utilisée pour le couple moyen ainsi que la FEM est la fonction Matérn 5/2. Pour le ratio des puissances et vu que c'est irrégulier, on a utilisé la fonction exponentielle. Le plan d'expériences initial est un Hypercube Latin optimisé avec un critère « Maximin » car il présente de bonnes propriétés de projection (sur chacun des axes d'entrée) et est optimisé pour couvrir au mieux l'espace de recherche [7].

### 3.2. Méthodologies d'optimisation

Pour l'optimisation basée sur des modèles de substitution fixes (Optim1), nous avons utilisé l'algorithme génétique NSGA II [8] qui a été couplé aux métamodèles créés.

Pour Optim2 et Optim3, une approche BO a été adoptée. Cette approche est souvent dédiée à l'optimisation des fonctions boîte noire dont les évaluations sont coûteuses en termes de calcul. Elle s'appuie sur des GPs construits à partir de données simulées $D_n = \{x_i, f(x_i)\}_{i=1}^n$. Comme décrit précédemment, ces GPs sont notés $\hat{f}_1$, $\hat{f}_2$, et $\hat{g}$. Ensuite, de nouveaux points de simulation $x_{new} = \{x_{new}^j\}_{j=1}^q$ sont choisis en fonction d'un critère de remplissage (infill criterion) basé sur les modèles GP courants, plus précisément sur leurs moyennes et leurs variances. L'objectif de ce critère est double : explorer le domaine de conception pour rendre le modèle GP plus prédictif dans les zones inexplorées (objectif d'exploration) et également dans les zones présentant un grand potentiel par rapport aux objectifs d'optimisation (exploitation). Ainsi, l'optimisation de ce critère de remplissage fournit à chaque itération des points à simuler. Les nouvelles simulations $Y_{new}$ sont ainsi utilisées pour mettre à jour les modèles GP. Cette procédure est répétée jusqu'à ce qu'un critère d'arrêt soit satisfait. Le critère d'arrêt est souvent défini comme un budget maximal de simulations. Cette approche est schématisée à la Fig. 2.

Dans cette étude, nous avons utilisé deux critères d'enrichissement : Pour Optim2, on a utilisé q-batch Pareto Efficient Global Optimization (qParEGO) [9]. Pour Optim3, c'était q-batch Expected Hypervolume Improvement (qEHVI) [9]. Nous décrivons brièvement ces méthodes dans les sections suivantes.

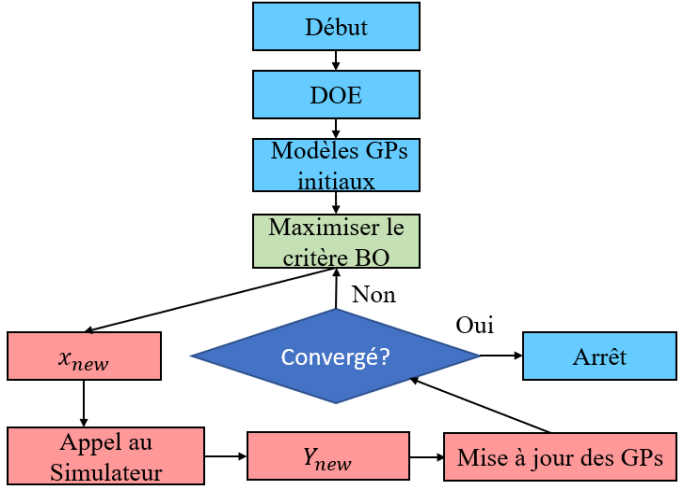

Fig. 2. Approche d'optimisation bayésienne.

*3.3. q-batch Pareto Efficient Global Optimization*

L'algorithme Pareto Efficient Global Optimization (ParEGO) [10] et sa version à q points (qParEGO) est une généralisation directe de l'algorithme bien connu Efficient Global Optimization (EGO) [11] aux problèmes d'optimisation multi-objectifs. L'algorithme EGO vise à résoudre le problème suivant :

$$\min_{\mathbf{x} \in D_x} f(\mathbf{x}) \quad (3)$$

L'algorithme EGO est illustré par la Fig. 3. La fonction à minimiser, $f$, en pointillés noirs n'est connue que par un nombre limité de simulations représentées par les points rouges. A partir de ces points, un processus gaussien $\hat{f}$ est construit : sur la figure, sa moyenne est représentée par la courbe bleue et un intervalle de confiance de ce métamodèle (en gris) est déduit de l'écart type du GP ($\pm 1.96\sqrt{C(x,x)}$). La ligne pointillée rouge représente la valeur minimale/optimale trouvée avec les simulations (points rouges). Nous définissons ensuite un point $x$ en bleu, candidat pour une nouvelle simulation. En ce point, nous pouvons prédire, grâce au GP $\hat{f}$, une amélioration par rapport au point minimum actuel (*Predicted improvement for the candidate for a new simulation using current kriging model*). Mathématiquement, cette fonction est définie comme suit :

$$I(x) = \left( \min_{i \in 1,2,\ldots n}(f(x_i)) - \hat{f}(x) \right)^+ \quad (4)$$

avec $(x)^+ = \max(x, 0)$. Ainsi, si la prédiction d'un point quelconque n'améliore pas le minimum actuel ($\hat{f}(x) \geq \min_{i \in 1,2,\ldots n}(f(x_i))$), alors son amélioration est égale à zéro.

Nous ne connaissons pas la valeur réelle de l'amélioration obtenue avec $x$ puisque nous n'avons pas évalué f au point $x$ à ce stade, mais nous pouvons estimer une moyenne de l'amélioration prédite par le GP à partir de la densité de probabilité associée au GP au point $x$, qui est la zone verte sous la courbe gaussienne autour de $x$. Nous définissons donc l'amélioration espérée (EI : *Expected Improvement*) comme l'espérance d'améliorer la meilleure valeur simulée courante :

$$EI(x) = \mathbb{E}_{\hat{f}}[I(x)|D_n] = \int_{-\infty}^{\infty} I(x) \, d\mathbb{P}_{\hat{f}} \quad (5)$$

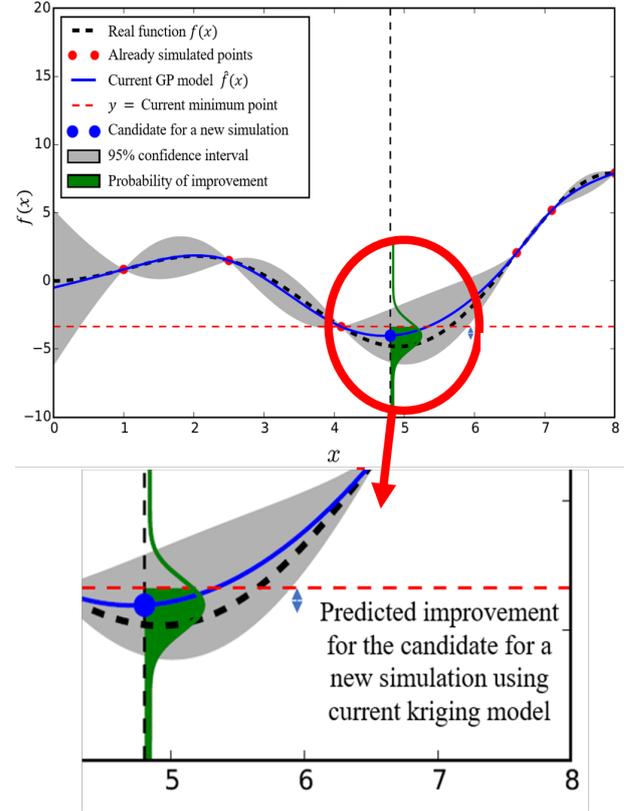

Fig. 3. Illustration de la méthode EGO pour un problème de minimisation. Figure adaptée de [12].

où $\mathbb{E}_{\hat{f}}$ est l'espérance par rapport à la densité de probabilité déduite du GP associée à $\hat{f}$ évaluée en x, et $d\mathbb{P}_{\hat{f}}$ est la mesure de cette densité. Cette expression peut être calculée analytiquement [11].

Dans notre étude, nous traitons un problème multi-objectifs. La première approche pour résoudre un problème multi-objectifs consiste à résoudre plusieurs problèmes mono-objectifs, par exemple, définis avec une nouvelle fonction objectif F(x) représentant une combinaison (généralement une somme pondérée) des fonctions objectives originales. Plus précisément, ParEGO utilise la fonction suivante :

$$F(x) = \alpha\big(wf_1(x) + (1-w)f_2(x)\big) + \max(wf_1(x), (1-w)f_2(x)) \quad (6)$$

avec $\alpha$ généralement égal à 0,05 [10]. Grâce au terme non-linéaire $\max(wf_1(x), (1-w)f_2(x))$, nous pouvons obtenir des solutions dans les zones non convexes du front Pareto de $f_1$ et $f_2$ [10]. Ces zones ne seraient pas accessibles avec une optimisation uniquement de $\alpha(wf_1(x) + (1-w)f_2(x))$. Quant au poids w, ses valeurs sont choisies suivant une loi uniforme entre 0 et 1 et varient à chaque itération de l'algorithme afin de mieux couvrir le front de Pareto. Une fois que F(x) est définie, nous pouvons utiliser le critère EI sur celle-ci. Le critère EI peut être étendu à un critère « multipoint » : il s'agit d'ajouter (q>1) points à chaque itération au lieu d'un seul ce qui est important si on l'on vise à réaliser des simulations en parallèle. Cela peut être fait en utilisant une stratégie *greedy* : optimiser le critère EI pour le premier point, conditionner les GPs sur le point déjà trouvé (en utilisant une observation virtuelle obtenue grâce à un échantillonnage de réalisations des GPs) et optimiser à nouveau le critère EI pour trouver le point suivant. Cette technique est répétée jusqu'à ce que $q$ points soient sélectionnés.

### 3.4. *q-step Expected Hypervolume Improvement*

Une autre approche d'optimisation Bayésienne dédiée à l'optimisation multi-objectifs est basée sur l'amélioration espérée de l'hypervolume (EHVI : Expected HyperVolume Improvement). Par analogie avec le critère EI, un critère d'amélioration en espérance est maximisé afin de choisir le point $x$ pour la prochaine simulation. Pour la méthode EHVI, la "qualité" de la solution actuelle (donnée par l'ensemble de Pareto actuel défini par les points noirs sur la Fig. 4, noté $\text{Par}(\{x_i\}_{i=1}^n)$) est mesurée par l'hypervolume dominé (surface en bleu sur la Fig. 4, notée $\text{HV}(\text{Par}(\{x_i\}_{i=1}^n))$) par rapport à un point de référence $f_{HV}$ (point rouge sur la Fig. 4). Plus l'hypervolume dominé est grand, meilleure est la solution. L'amélioration est alors obtenue en calculant l'hypervolume "gagné" (surface en vert dans la Fig. 4) lors de l'ajout d'un nouveau point (point vert dans la Fig. 4). Nous définissons l'amélioration de l'hypervolume comme étant :

$$\text{HVI}(x) = \text{HV}\left(\text{Par}(\{x_i\}_{i=1}^n) \cup \hat{f}_J(x)\right) - \text{HV}\left(\text{Par}(\{x_i\}_{i=1}^n)\right) \quad (7)$$

où nous notons $\hat{f}_J$ le modèle GP joint pour les 2 objectifs. Comme dans le cas de l'EI, nous ne connaissons pas la valeur réelle de l'hypervolume vert au point candidat x car nous n'avons pas calculé $f_1$ and $f_2$ à ce point. Nous pouvons néanmoins estimer l'espérance de cette amélioration de l'hypervolume grâce aux GPs $\hat{f}_1$ et $\hat{f}_2$ associés à $f_1$ et $f_2$:

$$\text{EHVI}(x) = \mathbb{E}_{\hat{f}_J}[\text{HVI}(x)] = \int_{-\infty}^{\infty} \text{HV I}(x) d\mathbb{P}_{\hat{f}_J} \quad (8)$$

Comme suggéré précédemment, le critère de remplissage EHVI peut être converti en un critère multipoint grâce à une stratégie *greedy*.

### 3.5. *Gestion des contraintes*

La manière la plus connue de traiter les contraintes $g(x) \leq 0$ en optimisation Bayésienne est d'utiliser l'amélioration faisable (FI : feasibility improvement) :

$$\text{FI}(x) = I(x) 1_{[g(x) \leq 0]} \quad (9)$$

où $1_{[g(x) \leq 0]}$ est la fonction indicatrice sur $g(x) \leq 0$. Elle est égale à 1 si $g(x) \leq 0$ et à 0 sinon. Ce faisant, le critère d'EI devient l'amélioration faisable en espérance (EFI) :

$$\text{EFI}(x) = \mathbb{E}\big[I(x) 1_{[g(x) \leq 0]}\big] = \text{EI}(x) \, \mathbb{P}_{\hat{g}}(g(x) \leq 0) \quad (10)$$

où $\mathbb{P}_{\hat{g}}(g(x) \leq 0)$ est la probabilité de satisfaire les contraintes, obtenue à partir du GP $\hat{g}$ associé à g. Des expressions similaires peuvent être facilement obtenues avec EHVI. Les critères qParEGO et qEHVI utilisés dans ce travail sont calculés avec Botorch, un paquet codé en Python [14]. Dans Botorch, la fonction indicatrice est remplacée par une fonction sigmoïde pour régulariser [9].

### 3.6. *Processus d'optimisation*

Les méta-modèles fixes pour Optim1 ont été construits en utilisant un plan d'expériences LHS maximin de 450 points. Pour les 2 optimisations Optim2 et Optim3, nous avons utilisé la même approche LHS maximin mais avec 250 points. Une fois le DOE initial de 250 points obtenu, nous avons effectué 50 itérations de BO ; à chaque itération, nous avons optimisé les deux critères BO en ajoutant 4 points, ce qui a permis d'ajouter un total de 200 points pour chaque BO. D'un point de vue temps de calcul, les 3 approches font le même nombre d'appels (450) du simulateur EF ce qui permet une comparaison équitable.

## 4. RESULTATS DES OPTIMISATIONS

La Fig. 5 montre les points ajoutés $x_{new}$ tout au long des itérations pour les optimisations bayésiennes qParEGO et qEHVI.

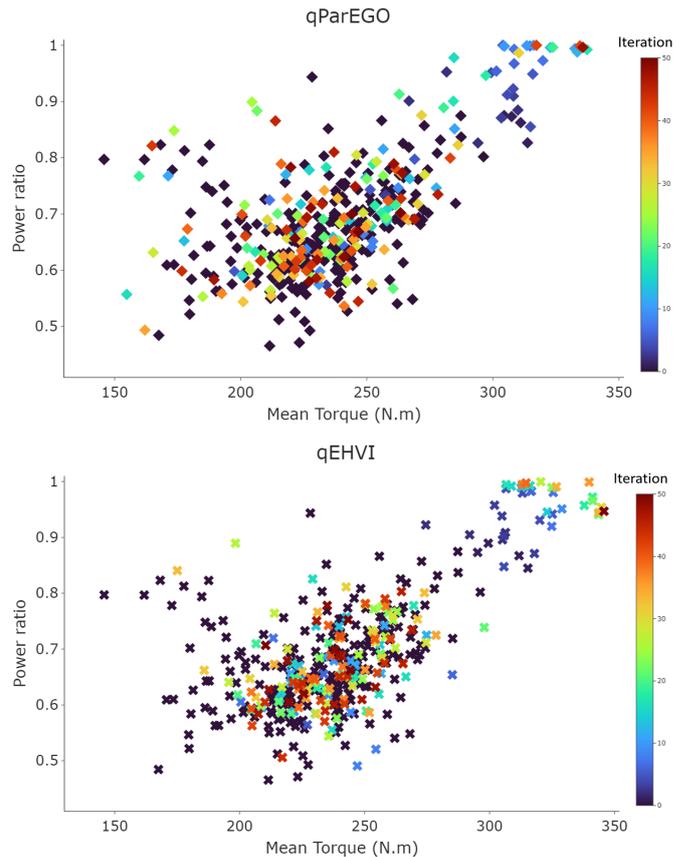

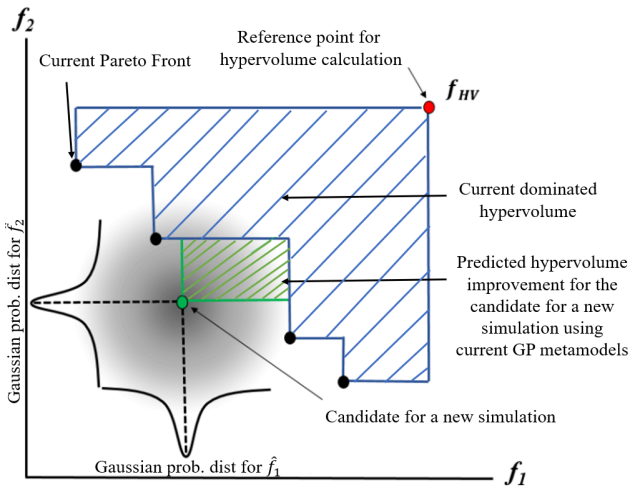

Fig. 4. Représentation de l'EHVI pour un problème de minimisation bi-objectifs. Figure adaptée de [13].

Fig. 5. Points de simulation ajoutés au fil des itérations (l'itération zéro correspond au plan d'expériences initial). Le graphique supérieur correspond à l'algorithme qEHVI et le graphique inférieur à l'algorithme qParEGO.

Cette figure montre que les deux approches BO ont trouvé plusieurs points dans la zone d'intérêt (coin supérieur droit), ainsi que des points repartis dans l'espace couple/ratio de puissance, ce qui montre les propriétés d'exploration et exploitation des critères BO.

La Fig. 6 compare les résultats obtenus avec les 3 différentes approches d'optimisation mises en œuvre dans ce travail. Les points en bleu foncé représentent le front de Pareto de Optim1, l'approche avec des modèles de substitution fixes. Afin de valider ces résultats prédits, on a choisi quelques points sur ce front de Pareto (cyan) pour recalculer leurs performances via des simulations par éléments finis (points en orange). Comme prévu, nous remarquons des différences entre les performances prédites et celles simulées. En effet, les modèles de substitution ont une précision limitée par un apprentissage sur un unique plan d'expériences choisi sur un critère d'exploration. Pour expliquer l'écart entre les performances prédites et celles simulées pour Optim1, on montre sur la Fig. 7 le front Pareto prédit ainsi que les 450 points du DOE. Il est clair que les performances des machines du front de Pareto sont une extrapolation des performances simulées des machines du plan d'expérience.

La Fig. 6 montre aussi les machines non dominées trouvées avec les approches qEHVI et qParEGO (points rouges et noirs, respectivement). Ces machines montrent de meilleures performances que celles trouvées avec Optim1.

La Fig. 8 compare la géométrie de trois machines sélectionnées (et encerclées) sur les fronts de Pareto de la Fig. 6 : Machine a) pour 'Optim1 (326,2 N.m ; 0,98), Machine b) pour Optim2 (335,7 N.m ; 0,99) et Machine c) pour Optim3 (345,8 N.m ; 0,95). Chaque machine présente un « bon compromis » sur son front de Pareto. Le Tableau 3 présente les valeurs numériques des paramètres d'entrée ainsi que les performances des machines sélectionnées.

La machine c) présente le plus de couple avec 345,8 N.m mais possède aussi la plus grande FEM avec 650 Vpp. On remarque que la valeur maximale du couple augmente avec la valeur de la FEM. En analysant les paramètres géométriques des machines, on remarque une corrélation entre l'angle d'ouverture de la barrière Beta_L1 et le couple (ainsi que la FEM) : Un angle d'ouverture plus faible signifie un couple plus élevé.

Le rapport de puissance de la machine c) est égal à 0,95. Pour augmenter ce rapport, un couple plus faible est à accepter ce qui montre que les deux objectifs choisis sont contradictoires.

La puissance maximale de 182,7 kW est atteinte pour la machine b). Cette machine possède aussi la quantité d'aimants la plus élevée avec 1,84 kg.

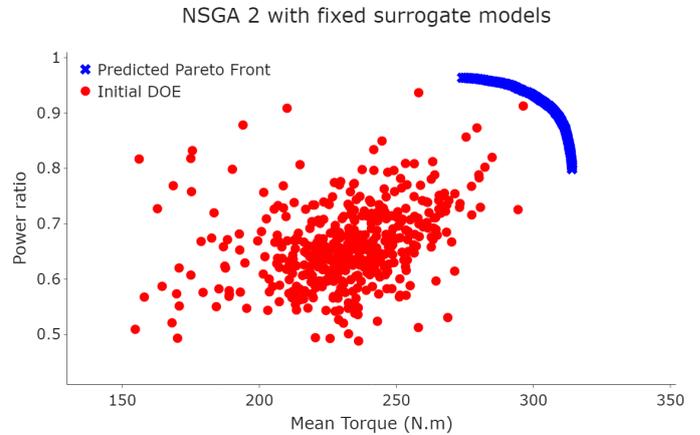

Fig. 7. Résultats des optimisations : Optim1 (points bleus), et le DOE initial associé (points rouges).

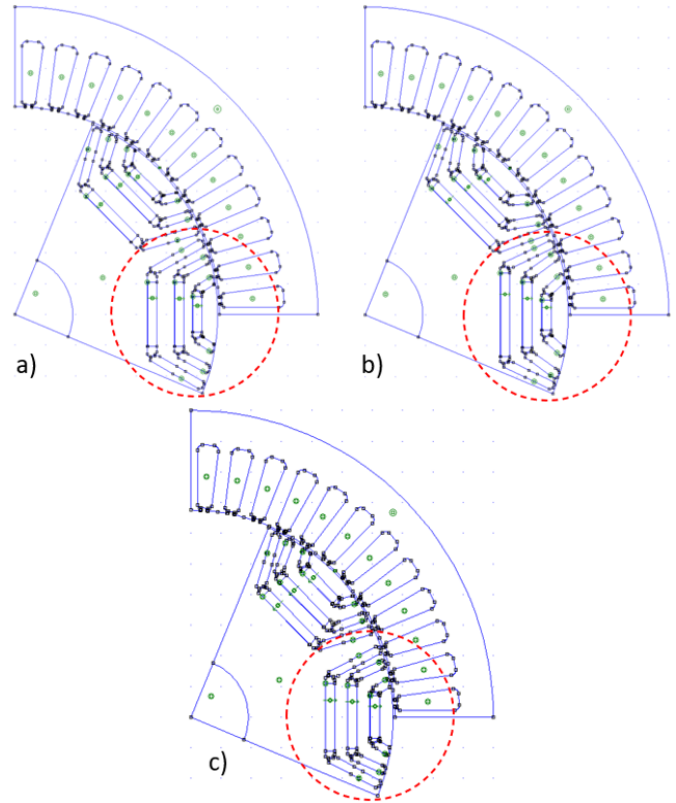

Fig. 8. Designs des machines pour des solutions sélectionnées dans la Fig. 6 : a) Optm1, b) Optim2, et c) Optim3.

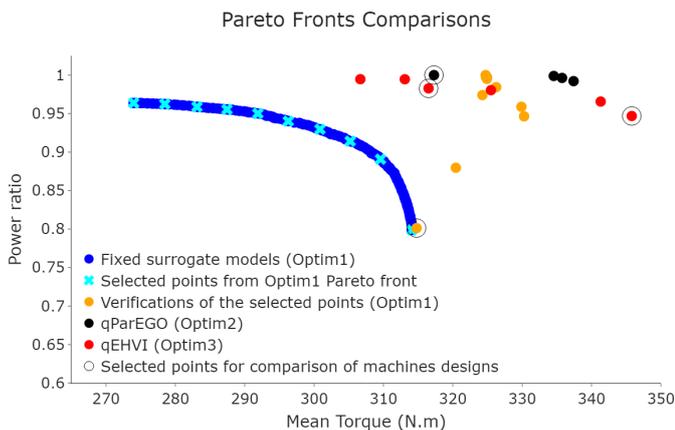

Fig. 6. Comparaison des résultats obtenus par les 3 méthodologies d'optimisation : avec modèles de substitution fixes : points bleus et cyans (points orange correspondent aux vérification EF des points cyans) et points non-dominées des stratégies bayésiennes qParEGO et qEHVI : points noirs et rouges respectivement.

Tableau 3: Résultats des optimisation : Variables de conception et quantités d'intêret.

|  |  | Machine a) | Machine b) | Machine c) |
|---|---|---|---|---|
| Variables de conception | Rad_PM_L1 | 0,86 | 0,89 | 0,89 |
|  | Rad_PM_L3 | 0,82 | 0,66 | 0,66 |
|  | Rad_PM_L3 | 0,64 | 0,6 | 0,8 |
|  | Rad_Brid_L1 | 0,2 | 0,1 | 0,24 |
|  | Rad_Brid_L2 | 0,17 | 0,23 | 0,1 |
|  | Rad_Brid_L3 | 0,43 | 0,51 | 0,38 |
|  | Beta_L1 | 0,87 | 0,83 | 0,67 |
|  | Beta_L2 | 0,41 | 0,72 | 0,29 |
|  | Beta_L3 | 0,53 | 0,7 | 0,5 |
|  | PM_Len_L1 | 0,78 | 0,9 | 0,86 |
|  | PM_Len_L2 | 0,84 | 0,9 | 0,85 |
|  | PM_Len_L3 | 0,77 | 0,7 | 0,9 |
| Quantités d'intêret | Couple | 326,2 N.m | 335,7 N.m | 345,8 N.m |
|  | Rapport de puissance | 0,98 | 0,99 | 0,95 |
|  | Puissance maximale | 180,2 kW | 182,7 kW | 177,4 kW |
|  | FEM à 8000 rpm | 624,6 Vpp | 629,9 Vpp | 649,6 Vpp |
|  | Poids aimants | 1,69 kg | 1,84 kg | 1,79 kg |

## 5. CONCLUSIONS

Nous avons présenté dans cet article une comparaison entre une optimisation par métamodèles fixes et deux approches bayésiennes appliqués au cas du dimensionnement du rotor d'une machine synchro-réluctante assistée par aimants permanents. Pour la première approche, nous avons remarqué des différences entre les performances des solutions prédites et simulées. Les approches bayésiennes basées sur des simulations éléments finis à chaque itération évitent cet écueil.

Grâce à la mise à jour adaptative des modèles de substitution, les approches bayésiennes permettent de trouver des solutions avec de meilleures performances, tout en utilisant un nombre limité de simulations éléments finis. Ce type d'approches peut être utile lorsqu'on est face à un problème d'optimisation nécessitant un nombre important de simulations ou bien utilisant un simulateur gourmand en temps de calcul comme dans des problèmes 3D, des problèmes multiphysiques ou bien des problèmes tenant compte des incertitudes.

## 6. REMERCIEMENTS

Nous voudrions remercier le CONAHCYT (Consejo Nacional de Humanidades Ciencias y Tecnologías) pour son support.

## 7. REFERENCES

[1] Ma, Cong; Qu, Liyan (2015) Multiobjective Optimization of Switched Reluctance Motors Based on Design of Experiments and Particle Swarm Optimization. In : IEEE Transactions on Energy Conversion, vol. 30, n° 3, p. 1144–1153. DOI: 10.1109/TEC.2015.2411677.

[2] Lei, Gang; Wang, Tianshi; Zhu, Jianguo; Guo, Youguang; Wang, Shuhong (2015) System-Level Design Optimization Method for Electrical Drive Systems—Robust Approach. In : IEEE Transactions on Industrial Electronics, vol. 62, n° 8, p. 4702–4713. DOI: 10.1109/TIE.2015.2404305.

[3] Reyes, Adan; Nasr, André; Sinoquet, Delphine; Hlioui, Sami "Robust design optimization taking into account manufacturing uncertainties of a permanent magnet assisted synchronous reluctance motor," 2022 IEEE International Conference on Electrical Sciences and Technologies in Maghreb (CISTEM), Tunis, Tunisia, 2022, pp. 1-6, doi: 10.1109/CISTEM55808.2022.10043885

[4] Zhang, Shen; Li, Sufei; Harley, Ronald G.; Habetler, Thomas G. (2018) An Efficient Multi-Objective Bayesian Optimization Approach for the Automated Analytical Design of Switched Reluctance Machines: 2018 IEEE Energy Conversion Congress and Exposition. 23-27 September 2018, Portland, OR, USA. 2018 IEEE Energy Conversion Congress and Exposition (ECCE). Portland, OR, USA, Piscataway, New Jersey: Institute of Electrical and Electronics Engineers, p. 4290–4295.

[5] Lei, Gang; Zhu, Jianguo; Guo, Youguang; Liu, Chengcheng; Ma, Bo (2017) A Review of Design Optimization Methods for Electrical Machines. In : Energies, vol. 10, n° 12, p. 1962. DOI: 10.3390/en10121962.

[6] C. E. Rasmussen and C. K. Williams, Gaussian Processes for Machine Learning. Cambridge: The MIT Press, 2006.

[7] Pronzato, Luc; Müller, Werner G. "Design of computer experiments: space filling and beyond." In : Statistics and Computing, vol. 22, n° 3, p. 681–701. 2012

[8] Kalyanmoy Deb, Amrit Pratap, Sameer Agarwal, and T. Meyarivan, "A Fast and Elitist Multiobjective Genetic Algorithm: A Fast and Elitist Multiobjective Genetic Algorithm: NSGA-II," *IEEE Transactions on Evolutionary Computation*, pp. 182–197. 2002.

[9] Daulton, Samuel; Balandat, Maximilian; Bakshy, Eytan (2020) Differentiable Expected Hypervolume Improvement for Parallel Multi-Objective Bayesian Optimization. In : Advances in Neural Information Processing Systems 33.

[10] J. Knowles. Parego: a hybrid algorithm with on-line landscape approximation for expensive multiobjective optimization problems. IEEE Transactions on Evolutionary Computation, 10(1):50–66, 2006.

[11] Donald R. Jones, Matthias Schonlau, William J. Efficient global optimization of expensive black-box functions. Journal of Global Optimization, 13, 1998.

[12] Guerra, Jonathan Optimisation multi-objectif sous incertitudes de phénomènes de thermique transitoire. PhD Thesis UNIVERSITEE DE TOULOUSE 2016.

[13] Palar, Pramudita S.; Zuhal, Lavi R.; Chugh, Tinkle; Rahat, Alma (01062020) On the Impact of Covariance Functions in Multi-Objective Bayesian Optimization for Engineering Design: AIAA Scitech 2020 Forum. AIAA Scitech 2020 Forum. Orlando, FL. Reston, Virginia: American Institute of Aeronautics and Astronautics.

[14] Balandat, Maximilian; Karrer, Brian; Jiang, Daniel R.; Daulton, Samuel; Letham, Benjamin; Wilson, Andrew Gordon; Bakshy, Eytan (2020) BoTorch: A Framework for Efficient Monte-Carlo Bayesian Optimization. In : Advances in Neural Information Processing Systems 33.